

\documentstyle{amsppt}
\magnification=\magstep 1
\TagsOnRight
\NoBlackBoxes
\rightheadtext{Universal arrows to forgetful functors}
\def\norm #1{{\left\Vert\,#1\,\right\Vert}}
\def\g {{\frak g}}
\def\Q {{\Bbb Q}}
\def\R {{\Bbb R}}

\def\N{{\Bbb N}}

\def\Z {{\Bbb Z}}

\topmatter

\title
Universal arrows to forgetful functors from
categories of topological algebra
\endtitle

\author Vladimir G. Pestov
\endauthor

\affil
Department of Mathematics\\
Victoria University of Wellington \\
P.O. Box 600 \\
Wellington, New Zealand
\\ \\
vladimir.pestov$\@$vuw.ac.nz
\endaffil
\abstract{We survey the present trends in
theory of universal arrows to forgetful functors
from various categories of
topological algebra
and functional analysis
to categories of topology and topological algebra.
Among them are free topological groups, free locally convex spaces,
free Banach-Lie algebras, and more. An accent is put on relationship of
those constructions with other areas of mathematics and their possible
applications. A number of open problems
is discussed; some of them belong to universal arrow theory,
and other
may become amenable to the methods of this theory.}
\endabstract
\subjclass{18A40, 22A05, 15A75, 16B50, 16W50, 17B01, 17B35,
22A30, 22E65, 43A40, 46A03, 46B04, 46H20, 46L05, 46M15,
54H11, 54A50, 54B30}
\endsubjclass
\endtopmatter
\document

\item{}{\sl  Introduction   (2)}
\item{1.}{\sl  Major classical examples  (3)}
\item{2.}{\sl  Structure of free topological groups  (7)}
\item{3.}{\sl  $M$-equivalence and dimension  (11)}
\item{4.}{\sl  Applications to general topological groups  (15)}
\item{5.}{\sl  Free products of topological groups  (17)}
\item{6.}{\sl  Free Banach-Lie algebras and their Lie groups  (19)}
\item{7.}{\sl  Lie-Cartan theorem  (20)}
\item{8.}{\sl  Locally convex Lie algebras and groups  (21)}
\item{9.}{\sl  Supermathematics  (23)}
\item{10.}{\sl $C^\star$-algebras and noncommutative mathematics  (25)}
\item{}{\sl Acknowledgments   (26)}
\item{}{\sl Bibliography  (28)}
\pagebreak

\heading
Introduction
\endheading
\smallpagebreak

The concept of a universal arrow was invented by P. Samuel in 1948 \cite{Sa}
and put in connection with his investigations on free topological groups.
The following definition is taken from the book \cite{MaL}.

\proclaim{Definition}
If $S:D\to S$ is a functor
and $~c~$ an object of $C$, a universal arrow from $c$ to $S$ is a pair $<r,u>$
consisting of an object $r$ of $D$ and an arrow $u: c\to Sr$ of $C$,
such that to every pair $<d,f>$ with $d$ an object of $D$
and $f:c\to Sd$ an arrow of $C$, there is a
unique arrow $f':r\to d$ of $D$ with
$Sf'\circ u =f$.

In other words, every arrow $f$ to $S$ factors uniquely through the
universal arrow $u$, as in the commutative diagram

$$
\CD
c @>u>> Sr \\
@| @VV Sf' V \\
c @>f>> Sd \\
\endCD
$$

\qed\endproclaim

This notion bears enormous generality and strength, and
at present it is
an essential ingredient of category theory \cite{MaL}
and theory of toposes \cite{Joh}.
In fact, many mathematical constructions can be
interpreted as universal arrows of one of another kind.
Examples are: quotient structures and substructures,
products and coproducts, including
algebraic and topological tensor products,
universal enveloping algebras, transition from a Lie algebra to a
simply connected Lie group and {\it vice versa}, compactifications
of all kinds (Stone-\v Cech, Bohr, and others), completions,
prime spectra of rings,
and much more.

We are interested in the particular case where
$S$ is a {\it forgetful functor}
from some category of topological algebra or functional analysis, $D$,
to another category of topological algebra
or functional analysis or a category of topology, $C$.
Historically the first, and studied in most detail, is the
construction of the {\it free topological group}, $F(X)$, over a topological
space $X$, where  $C$ is the category of Tychonoff topological spaces
and continuous mappings
and $D$ is
the category of Hausdorff topological groups and continuous homomorphisms.
A number of similar constructions have received a
comprehensive treatment, among them are
free Abelian topological groups, free compact groups,
free locally convex spaces.
At the same time, in recent years similar constructions have arisen
--- either explicitely or implicitely ---
in other areas of mathematics. In some cases
no attempt has been made to establish a bridge between those and former
types of universal arrows --- although seemingly such a connection would
facilitate a study of new constructions.
Among the disciplines where
new types of universal arrows to forgetful
functors are likely to play a noticeable role,
are infinite-dimensional Lie theory,
supermanifold theory, differential geometry,
$C^\star$-algebras and ``quantized'' functional analysis.

We do not aim at a comprehensive presentation of the subject
outlined in the title of this paper, nor we give detailed
proofs of the results: such an elaborate approach would lead to
a voluminous treatise.
Instead, we discuss a few carefully selected lines of
development which, as we see it, dominated the research over
more than 50 years. We are focussing on the most interesting
unsolved problems.
Also, we do our best in forecasting the future directions of the theory,
paying special attention to recent germs of it
in areas of mathematics
bordering topological algebra (Lie theory,
functional analysis and mathematical
physics).

This small survey inevitably tends to the results and ideas
coming from the Russian (or, in a more politically correst
language, ex-Socialist,
to cover Ukrainian, Moldavian, Georgian, Bulgarian
and other contributors) school of
universal arrow theorists, where the author himself comes from.
Most probably and to the author's regret, the
contributions from the other two major centers ---
the Australian and the American schools --- were
underrepresented in this article. As a matter of fact,
the author's personal tastes
and research work of his own were prevalent in selecting topics
for discussion.

Our bibliography, although (intentionally) not complete,
is hopefully ``everywhere dense'' in the subject
(a comparison due to Kac \cite{Kac1}).

\smallpagebreak
\heading
1. Major classical examples
\endheading
\smallpagebreak

The following are major
examples of universal arrows to forgetful functors
from categories of topological algebra, which are subject of
a traditional study in this area.
We are marking with a lozenge ($^\lozenge$)
those notions which will be later considered in this survey to some
extent.
By abuse of terminology and notation, we will
sometimes identify a universal arrow with its target object
(no confusion should however result from that).

\smallskip
\item {\bf 1.$^\lozenge$}
$C=\bold{Tych}$ (the category of Tychonoff topological spaces
and continuous mappings)
 and $D=\bold{TopGrp}$
(the category of Hausdorff topological groups and continuous homomorphisms).
The universal arrow from an object
$X\in C$ (a Tychonoff space) to the forgetful functor $S:D\to C$ is
the {\it (Markov) free topological group over $X$},  $F(X)$.
\smallskip
This notion was introduced in 1941 by Markov \cite{Mar1}  who presented his
results in most detail somewhat later \cite{Mar2}.
Among those mathematicians who responded first to the new concept, were
Nakayama \cite{Nak}, Kakutani \cite{Kak}, Samuel \cite{Sa}
and Graev \cite{Gr1};
the latter work  has had an
enormous impact on later investigations
in the area, and the paper by Samuel, as we have already mentioned,
has produced a deep methodological insight.
\smallskip
\item {\bf 2.~} $C=\bold{Tych}_\ast$
(the category of pointed Tychonoff topological spaces
and continuous mappings preserving base points)
and $D=\bold{TopGrp}$
(the base point of a topological group being $e$, the identity).
The universal arrow from an object
$X\in C$ (a pointed Tychonoff space) to the forgetful functor $S:D\to C$ is
the {\it Graev free topological group over $X$},  $F_G(X)$.
\smallskip
In fact, the Markov and Graev free topological groups are very closely
related to each other by means of the following
short exact sequence:

$$ e \to \Z \to F_M(X) \to F_G(X) \to e$$

The choice of a basepoint $\ast\in X$ does not affect the topological
group $F_G(X)$. The Markov free group of $X$ is isomorphic to the
Graev free group of the disjoint sum $X\oplus\{\ast\}$. \cite{Gr1,2}.
This is why we consider the Markov free topological groups only.
Anyway, the Graev approach seems more convenient in some other cases
such as free Banach spaces and free Banach-Lie algebras over metric spaces.
\smallskip
\item {\bf 3.~} $C=\bold{Met}_\ast$ (the category of pointed metric spaces)
 and $D=\bold{MetGrp}$ (the category of groups endowed with bi-invariant
metric).
The universal arrow from an object
$(X,\rho, \ast)\in C$ (a pointed metric space) to
the forgetful functor $S:D\to C$ is
the {\it free group over $X\backslash \{\ast\}$
endowed with the Graev metric $\bar\rho$}.
\smallskip
This concept is due to Graev \cite{Gr1,2}. The metrized group
$(F(X), \bar\rho)$ is of no particular interest by itself; it
deserves attention as an auxiliary device for studying the free
topological group $F(X)$. An amazing example of such kind is
the Arhangel'ski\u\i 's theorem from \cite{Arh4}. If one wants to consider
Graev matrics on a Markov free group then one should start with a
fixed metric $\rho$ on the set $X\oplus\{e\}$.
\smallskip
\item {\bf 4.~} $C=\bold{Tych}$
and $D=\Cal{V}$ is a {\it variety} of Hausdorff topological groups,
considered as a subcategory of $\bold{TopGrp}$.
The universal arrow from an object
$X\in C$ (a Tychonoff space) to the forgetful functor $S:D\to C$ is
the {\it  free topological group over $X$ in the variety ${\Cal V}$},
$F_{\Cal V}(X)$.
\smallskip
Varieties of topological groups can be understood in different
sense (cf. \cite{Mo1,2,10} and \cite{Pr2,3,PrS}).
It would not be clear what is the ``right'' notion until a
non-disputable version of the Birkhoff theorem for topological
groups is obtained (see, however, \cite{Ta}).
Anyway, all of the most important classes
of topological groups fit both definitions.
Examples of varieties are: the variety of SIN groups
(topological groups with
equivalent left and right uniformities) \cite{MoTh1};
that of topological groups with quasi-invariant basis \cite{Kats1}
($=$ $\aleph_0$-balanced groups in \cite{Arh5}); of totally bounded,
or precompact,
groups; of $\aleph_0$-bounded groups \cite{Gu, Arh5} etc.
There is a survey on free topological groups in varieties  \cite{Mo10}.
A free topological group, $F_{\Cal V}(X)$, in a variety $\Cal V$
is actually the composition of the universal arrow $F(X)$ and
the universal arrow from $F(X)$ to the natural embedding functor
$\Cal V\to\bold{TopGrp}$.
The notion of a free topological group relative to classes of topological
groups, considered by Comfort and van Mill \cite{ComvM},
can be redefined in terms of free topological groups in relevant varieties,
and the questions of existence of such
free topological groups completely reduces  to certain questions about
free topological groups in  varieties.

The following is the most important particular case.
\smallskip
\item {\bf 5.~} $C=\bold{Tych}$
 and $D=\bold{AbTopGrp}$ (the category of Abelian topological groups and
continuous homomorphisms).
The universal arrow from an object
$X\in C$ (a Tychonoff space) to the forgetful functor $S:D\to C$ is
the {\it (Markov) free Abelian topological group over $X$},  $A(X)$.
\smallskip
\item {\bf 6.~} $C=\bold{Tych}_\ast$
 and $D=\bold{AbTopGrp}$.
The universal arrow from an object
$X\in C$ (a Tychonoff space) to the forgetful functor $S:D\to C$ is
the {\it Graev free Abelian topological group over $X$},  $A_G(X)$.
\smallskip
Of course, $A(X)$ (resp. $A_G(X)$) is just the abelianization of $F(X)$
(resp., $F_G(X)$).
\smallskip
\item {\bf 7.~} $C=\bold{Tych}$
 and $D=\bold{CompGrp}$ (the category of compact topological groups and
continuous homomorphisms).
The universal arrow from an object
$X\in C$ (a Tychonoff space) to the forgetful functor $S:D\to C$ is
the {\it free compact group over $X$},  $F_C(X)$.
\smallskip
Remark that the free compact group $F_C(X)$ is nothing but the
Bohr compactification, $bF(X)$, of the free topological group, $F(X)$.
(The Bohr compactification, $bG$, of a topological group $G$
\cite{Mo9} is
the universal arrow from $G$ to the embedding functor
$\bold{CompGrp}\to\bold{TopGrp}$.)
We do not touch free compact groups in our survey, and refer the
reader to the series of papers by Hofmann and Morris
\cite{Hf,HfMo1-5}.
Also free compact groups may be viewed as completions
of {\it free precompact groups} (or, just the same,
free totally bounded groups), that is, free topological groups in the
correspondent variety. Free precompact groups have been studied
recently in
connection with some questions of dimension theory \cite{Sh}.

Of course, the notion of the free compact Abelian group over $X$
also makes sense, and the structure of such groups has been described
in detail ({\it loc. cit.}).
\smallskip
\item {\bf 8.$^\lozenge$} $C=\bold{Unif}$ is the category of uniform spaces
and $D=\bold{TopGrp}$. There are at least four ``natural''
forgetful functors from $D$ to $C$; our choice as $S$ is the functor
$S$ assigning to a topological group $G$ the {\it two-sided} uniform
structure on it; we will denote the resulting uniform space by $G_t$.
The universal arrow from an object $X\in \bold{ Unif}$
(a uniform space) to the functor $S$
is the {\it free topological group over $X$}, or
the {\it uniform free topological group}, $F(X)$.
\smallskip
This was an invention of Nakayama \cite{Nak}. Free topological groups over
uniform spaces later proved to be a most natural framework for analysing
some aspects of free topological groups, see \cite{Nu2}.
Free topological groups over uniform spaces provide a
straightforward generalization of free topological groups over
Tychonoff spaces, because
for a Tychonoff space $X$ the free topological grop over $X$ is
canonically isomorphic to the free topological group over the
finest uniform space associated to $X$.

\smallskip
\item {\bf 9.~} By replacing the category $\bold{Tych}$ by $\bold{Unif}$ in
the items 2,4,6 one comes to the obviously defined concepts of
a {\it (Graev) free (Abelian) topological group over a uniform space}.
\smallskip
\item {\bf 10.~} $C=\bold{Tych}$ (resp., $\bold{Unif}$) and
$D=\bold{LCS}$ (the category of locally convex spaces and continuous linear
operators). The universal arrow from an object $X\in C$ (a Tychonoff space;
resp., a uniform space) to the forgetful functor
$S: D\to C$ (which in the second case is also
defined unambiguously, unlike in item 8) is the
{\it free locally convex space over a topological (uniform) space} $X$,
and is denoted by $L(X)$.
\smallskip
This concept is also an invention of Markov \cite{Mar1}.
However, for some reason it received no immedeate attention from
mathematical community until the paper by Ra\u\i kov \cite{Rai2}.
The most important of later developments is due to Uspenski\u\i\ \cite{U2}.
A particular case of this construction --- the notion of a vector space
endowed with finest
locally convex topology --- is well known in functional analysis \cite{Sch};
it is actually the free locally convex space over a
discrete topological space $X$.
\smallskip
\item {\bf 11.~} As in item 4, one can consider universal arrows from
an object of $\bold{Tych}$ to the forgetful functor ${\Cal V}\to \bold{Tych}$
where ${\Cal V}$ is a {\it variety} of locally convex spaces in one or
another sense. We denote the resulting {\it free locally convex space
over $X$ in the variety} ${\Cal V}$ by $L_{\Cal V}(X)$.
\smallskip
We refer the reader to a very solid paper \cite{DmOS} by
Diestel, Morris and Saxon, and a survey \cite{Mo10}
by Morris. Other references include \cite{Ber}.
\smallskip
\item {\bf 12.$^\lozenge$} If ${\Cal V}$ is the variety of
locally convex spaces
with weak topology then the resulting {\it free locally convex space
with weak topology} over a Tychonoff space $X$ is denoted by
$L_p(X)$.
\smallskip
This concept seemingly was well known in functional analysis
for decades,
because the space $L_p(X)$ is the weak dual of the space $C_p(X)$
of continuous functions on $X$ in the topology of pointwise
(simple) convergence. See, e.g., \cite{Wh} and references therein.
\smallskip
\item {\bf 13.~} $C=\bold{Met}_\ast$ and $D=\bold{Ban}$
is the category of
{\it complete normed linear spaces and linear operators of norm}
$\leq 1$.
The universal arrow from an object $X\in\bold{Met}_\ast$ to the forgetful
functor $S:D\to C$ (the origin is a base point) is the
{\it free Banach space over a pointed metric space}, $B(X)$.
\smallskip
This object first appeared in the paper by Arens and Eells \cite{ArE};
see also \cite{Rai2; Pe9}.
However, it was considered by functional analysts
independently and at a different angle of view: the normed space
$B(X)$ is known as the {\it predual
of the space $Lip~(X)$ of Lipshitz
functions} on a pointed metric space $X$.
\smallskip
\item{\bf 14.~} $C=\bold{Tych}$ and
$D$ is the category of universal topological algebras
of a given signature $\Omega$. In this case the universal arrow from a
space $X$ to the forgetful functor $D\to C$ is the
{\it free universal topological algebra} over $X$.

Such algebras were first considered by Mal'cev \cite{Mal} and
others \cite{Ta, Pr2, PrS}.
We will not touch them in our survey.

\item{\bf 15.~}  $C=\bold{Tych}$ and $D$ is the category of topological
associative rings or associative algebras.
The resulting {\it free topological rings} and
{\it free topological algebras} have been also considered
by Arnautov, Mikhalev, Ursul and others \cite{AMV}.

Later in our survey we will consider also a number of less
traditional examples of universal arrows to forgetful
functors.
All of them are universal arrows to forgetful functors of one or
another kind. The following notion,
that of free product of topological groups,
at first seems not to fit into this scheme.

\item{\bf 16.$^\lozenge$} Let $C=\bold{TopGrp}\times \bold{TopGrp}$,
$D=\bold{TopGrp}$, and let $S$ be the diagonal functor
$\bold{TopGrp}\to \bold{TopGrp}\times \bold{TopGrp}$.
(That is, $S(G)=(G,G)$.)
The universal arrow from a pair $(G,H)$ of topological groups
to the functor $S$ is called the {\it free product of $G$ and $H$}
and denoted by $G\ast H$. In other terms,  $G\ast H$ is just the
coproduct of $G$ and $H$ in the category $\bold{TopGrp}$.

Anyway, it is well known that this notion
(belonging to Graev \cite{Gr3}) is of the same nature as
that of a free topological group,
those constructions share a number of common
properties and indeed, it can be (if necessary)
reshaped as a universal arrow to an appropriate
{\it forgetful} functor.
Let $C=\bold{TopGrp}\times \bold{TopGrp}$ be as above, and let
$D$ denote the category of all topological groups with two
fixed subgroups.
Then $G\ast H$ can be viewed as the universal arrow from
a pair $(G,H)$ to the forgetdul functor from $D$ to $C$ which
forgets the first group and sends a triple
$(F,G,H)$ to $(G,H)$.

\item{\bf 17.$^\lozenge$} In an obvious way, the concept of the free
product can be defined for arbitrary families of topological
groups, $\{G_\alpha : \alpha\in A\}$.
This product is denoted by $\ast_{\alpha\in A}G_\alpha$.

In all the aforementioned cases, similar methods, which are actually
of a categorial nature, are used to prove the existence, uniqueness and a
number of other properties of universal arrows.
We will summon those results as follows.

\proclaim{1.1. Theorem}
\item{(1)} In all cases 1--17 the universal arrow exists and is unique.
\item{(2)} In all cases apart from 4, the universal arrow
is an isomorphic embedding.
\item{(3)} In case 4, the universal arrow
is a homeomoprhic embedding if the variety $\Cal V$ contains at least
one non totally path-disconnected topological group.
\item{(4)} In all cases apart from  4 and 7,
the image of the iniversal arrow is topologically closed.
\qed\endproclaim

\medpagebreak
\heading  2. Structure of free topological groups
\endheading
\smallpagebreak

Among the first, and most vital, questions to be asked
about any universal arrow to
forgetful functor from a category of topological algebra,
is the question of description of the algebro-topological
structure of the target object of this arrow.
In some cases such a description poses no serious problems,
but for most (especially noncommutative) examples it is
rather challenging.
Since this question seems to be best
investigated for free topological groups,
we find it necessary --- and very instructive --- to
survey the state of affairs in this area.

\subheading{ 1. Description of topology
}
The topology of a free topological group $F(X)$
is rather complicated, and among the achievements of Graev \cite{Gr1,2}
was  a description of the topology of $F(X)$ in the case where $X$
is a compact space.
Later his description was transferred to the so-called
$k_\omega$-spaces by Mack, Morris and Ordman \cite{MaMoO},
which result has substantially widened the sphere of applicability
of the original description. We will give it in the strongest form.

Denote by $\tilde X=X\oplus -X \oplus \{e\}$
the disjoint sum of a Tychonoff space $X$, of its topological copy
$-X=\{-x :x\in X\}$, and a one-point space $\{e\}$.
For each $n=0,1,2,\dots$ there is an
obvisouly defined  canonical continuous mapping
$i_n:\tilde X^n\to F(X)$.
Denote by $F_n(X)$ the subspace of $F(C)$
image of $i_n$; it is closed.
A topological space $X$ is called a $k_\omega${\it -space} if
it can be represented as a union of countably many compact subsets
$X_n$ in such a way that the topology of $X$ is a {\it weak} topology
with respect to the cover $\{X_n:n\in\N\}$, that is,
a subset $A\subset X$ is closed iff so are al intersections
$A\cap X_n,~n\in\N$.
Not only every compact space is a $k_\omega$ space;
so is every countable $CW$-complex, every locally compact space with
countable base etc.

\proclaim{2.1. Theorem {\rm (Graev-Mack-Morris-Nickolas)}}
Let $X$ be a $k_\omega$-space. Then every mapping $i_n$ is quotient,
and a subset $A$ of $F(X)$ is closed if and only if
so are all intersections $A\cap F_n(X)$.
In particular, $F(X)$ is a $k_\omega$-space.
\qed\endproclaim

The above theorem does not admit any noticeable further generalization,
apart from some openly pathological cases, such as
the spaces $X$ where every $G_\delta$ set is open
(the author, unpublished, 1981).
In fact, it was shown in \cite{FOST} that the mapping $i_3$ is not quotient
even
for $X=\Q$. Answering two question raised in this paper, the author
had proved the following result \cite{Pe7}.

\proclaim{2.2. Theorem}
Let $X$ be a Tychonoff space. The mapping $i_2$ is
quotient if and only if $X$ is a strongly collectionwise normal space
(that is, every neighbourhood of the diagonal in $X\times X$ is
an element of the universal uniform structure of $X$).
\qed\endproclaim

However, the following property of the mappings $i_n$
proved to be extremely useful.

\proclaim{2.3. Theorem {\rm (Arhangel'ski\u\i\ \cite{Arh2,3})}}
Let $Y$ be any subset of $\tilde X^n$ such that
$i_n^{-1}i_n (Y)=Y$.
Then $i_n\vert_Y$ is a homeomorphism.
\qed\endproclaim

A very substantial body of results concerning the structure
of free topological groups over $k_\omega$ spaces have been deduced
(mostly by Australian and American
mathematicians) from Theorem 2.1
\cite{Br, BrH, F, FRST2, HMo, Katz, KatzMo1,2, KatzMoN1-4, Mo10,11, Nic1,3,
OrST}.

The following charming
theorem of Zarichny\u\i\ \cite{Zar1,2} puts
free topological groups in connection with infinite-dimensional
topology. The original result was stated for free Graev topological
groups, but it extends to free Markov groups immedeately
because topologically
the group $F(X)$ is a disjoint sum of countably many copies of
$F_G(X)$.

\proclaim{2.4. Theorem {\rm (Zarichny\u\i\ \cite{Zar1})}}
Let $X$ be a compact absolute neighbourhood retract and
$0<dim~X<\infty$.
Then the free topological group $F(X)$ and the
free Abelian topological group $A(X)$ are homeomorphic to an
open subset of the locally convex space with finest
topology $\R^\omega=\varinjlim \R^n$.
\qed\endproclaim

Returning back to general Tychonoff spaces
$X$, one can still describe the topology of $F(X)$ with the help of
mappings $i_n$,
but in a rather non-constructuve way.
The following construction have been performed by Mal'cev
\cite{Mal}.
Denote by $\frak T_0$ the quotient topology on $F(X)$
with respect to the direct sum of the mappings $i_n,~n\in\N$
from the space $\oplus_{n\in\N}\tilde X^n$.
It is Hausdorff but not necessarily a group topology.
Now construct recursively a transfinite chain of topologies
$\frak T_\lambda$ on $F(X)$ by defining $\frak T_{\lambda +1}$
as the quotient of the  topology
on $F(X)\times F(X)$ with respect to the
 mapping $(x,y)\mapsto x^{-1}y$,
and $\frak T_{\tau}$ for a limit cardinal $\tau$ as
the infimum of the chain of topologies $\frak T_{\lambda },~\lambda<\tau$.
It is clear that for some $\lambda$ large enough, the topology
$\frak T_{\lambda }$ coincides with the topology of $F(X)$.
Denote the least $\lambda$ with this property by $\lambda (X)$.
The following question is open for more than 30 years.

\proclaim{Problem {\rm (Mal'cev \cite{Mal})}}
Which values can $\lambda (X)$ assume?
\endproclaim

Seemingly, all one knows is that $\lambda (X)=1$ for $k_\omega$-spaces,
and  $\lambda (X)> 1$ for most spaces beyond this class
(for instance, for $X=\Q$).

Another long-standing problem asked by Mal'cev in the same paper \cite{Mal}
--- that of finding a constructive description
of a neighbourhood system of identity of a free topological group ---
had been solved by
Tkachenko \cite{Tk4}.
Later simpler versions of the Tkachenko's theorem have been obtained
by the author \cite{Pe7} and Sipacheva \cite{Si1}.
We will give one of the possible forms of the result.
It is more reasonable to put it for free topological groups $F(X)$
over {\it uniform spaces} (bearing in mind that
for a Tychonoff space $X$ the free topological grop over $X$ is
canonically isomorphic to the free topological group over the
universal uniform space associated to $X$).
Let $X=(X,\Cal U_X)$ be a uniform space.
Denote by $j_2$ a mapping from $X^2$ to $F_2(X)$ of the form
$(x,y)\mapsto x^{-1}y$, and by ${j^\ast}_2$ ---
a similar mapping of the form
$(x,y)\mapsto xy^{-1}$. If $\Psi\in (\Cal U_X)^{F(X)}$
is a family of entourages of diagonal indexed by elements of
the free group over $X$, then we put

$$\Cal V_\Psi=_{def} \cup\{ x\cdot [j_2(\Psi (x))\cup
{j^\ast}_2(\Psi (x))]\cdot x^{-1} : x\in F(X)\} $$

\noindent If $B_n$ is a sequence of subsets of some group then,
  following \cite{RoeD}, we denote

$$[(B_n)]=_{def} \cup_{n\in\N}\cup_{\pi\in S_n}
B_{\pi (1)}\cdot B_{\pi (2)} \cdot \dots
\cdot B_{\pi (n)}, $$

\noindent where $S_n$ is a symmetric group.

\proclaim{2.5. Theorem {\rm (Pestov \cite{Pe7})}}
Let $(X, \Cal U_X)$ be a uniform space.
A base of neighbourhoods of identity in the free topological group
$F(X)$ is formed by all sets of the form
$[(\Cal V_{\Psi_n})]$,
where $\{\Psi_n\}$ runs over the family of all countable sequences
of elements of $(\Cal U_X)^{F(X)}$.
\qed\endproclaim

\subheading{
2. Free subgroups
}
If $X$ is a subset of a set $Y$, then the free  group over the set of
generators $X$ is a subgroup of the free group over $Y$.
Now let $X$ is a topological subspace of a Tychonoff space $Y$;
there is still a canonical continuous group monomorphism
$F(X)\hookrightarrow F(Y)$, but it need not be a topological
embedding. For the first time it was noticed by Hunt and Morris
\cite{HuM}, and the example was $X=(0,1),~Y=[0,1]$.
Earlier Graev has shown \cite{Gr1,2} that if $Y$ is compact and $X$ is
closed in $Y$ then $F(X)\hookrightarrow F(Y)$ is in isomorphic
embedding of topological groups. This result was
transferred to $k_\omega$-spaces.
In \cite{Pe1,2,4} and \cite{Nu2} it was noticed independently
that a necessary condition for the monomorphism
$F(X)\hookrightarrow F(Y)$ to be topological is the
property that the restriction
$\Cal U_Y\vert_X$ of the universal uniformity $\Cal U_Y$ from $Y$ to $X$
coincides with the universal uniformity $\Cal U_X$ of $X$.
(It is just an immedeate consequence of the fact that
both left and two-sided uniformilies on $F(X)$
induce on $X$ its universal uniform structure --- the fact which in its
turn follows from existence of Graev's pseudometrics on $F(X)$
and was essentially known to Graev.)
In the same works \cite{Pe1,2,4} and \cite{Nu2} it was
shown,
answering a question by
Hardy,  Morris and Thompson
 \cite{HMoTh} that the above condition $\Cal U_Y\vert_X=\Cal U_X$
is sufficient in the case where $X$ is dense in $Y$.
A final positive answer was obtained by Uspenski\u\i\ \cite{U5}
after a series of results of intermediate strength \cite{U3}.

\proclaim{2.6. Theorem {\rm (Uspenski\u\i\ \cite{U5})}}
Let $X$ be a topological subspace of a Tychonoff space $Y$.
Then the monomorphism
$F(X)\hookrightarrow F(Y)$ is a topological embedding
if and only if $\Cal U_Y\vert_X=\Cal U_X$.
\qed\endproclaim

A different problem has been treated by Australian and American
universal arrow theorists for a long time.
Let $X$ and $Y$ be some particular topological spaces; in which cases
the free (Abelian) topological group over $X$ can be embedded
(not necessarily in a ``canonical'' way)
as a
topological subgroup into the free (Abelian) topological group over
$Y$? The main device under this approach was the above Theorem 2.1.
We will mention just one astonishing result in this direction.

\proclaim{2.7. Theorem {\rm (Katz and Morris \cite{KatzMo2})}}
If $X$ is a countable CW-complex of dimension $n$, then
the free Abelian topological group on $X$ is a closed subgroup of
the free Abelian topological group on the closed ball $B^n$.
\qed\endproclaim

\subheading{
3. Completeness
}
Our next topic can  be also traced back to Graev's papers \cite{Gr1,2}.
Graev has deduced from his description of topology of
the free group over a compact space that any such free topological
group is Weil complete
(that is, complete with respect to the left uniform structure).
The result remains true for free topological groups over
$k_\omega$-spaces.

Examples of  topological groups which are complete in their
two-sided uniformity but not Weil complete (and therefore admit
no Weil completion at all) are known for decades, but seemingly it
remains unclear whether free topological groups admit
Weil completion.
This question was asked by Hunt and Morris \cite{HuM}.
An obvious necessary condition for a free
topological group to be Weil-complete is the Dieudonn\'e completeness of
$X$, that is, completeness of $X$ w.r.t. the finest uniformity
$\Cal U_X$. The state of affairs with Weil completeness
is still unclear and one has only a series of partial results
stating the Weil completeness of free topological groups over
particular spaces \cite{U3}.

However, it seems in a sense more natural to examine
free topological groups for another form of completeness
--- the completeness w.r.t. the two-sided uniformity
(sometimes also called Ra\u\i kov completeness), \cite{Rai1}.
There exists a fascinating
comprehensive result for the completeness of this kind,
and the question about the validity of
such a result was first asked independently by Nummela \cite{Nu2}
and the author (in oral form, talk at the
Arhangel'ski\u\i 's seminar on topological
algebra at Moscow University,
February 1981).

\proclaim{2.8. Theorem {\rm (Sipacheva, \cite{Si2})}}
The free topological group $F(X)$ over a Tychonoff space $X$
is complete if and only if $X$ is Dieudonn\'e complete.
\qed\endproclaim

The idea of the proof is based on the notion of a special
universal arrow, $F_\rho (X)$, introduced by Tkachenko \cite{Tk3}.
Say that a subspace $Y$ of a topological group $G$ is
{\it Tkachenko thin} if for every neighbourhood of identity, $U$,
the set $\cap \{yUy^{-1}:y\in Y\}$
is a neighbourhood of identity.
Consider the category of pairs $(G,Y)$ where $G$ is a
Hausdorff topological group and $Y$ is a Tkachenko thin subset of $G$, and
obvious morphisms between them, and let $S$ be the functor
from this category to $\bold{Tych}$ of the form $(G,Y)\mapsto Y$.
Now by $F_\rho (X)$ one denotes
the universal arrow from a Tychonoff space $X$ to
the functor $S$. There is a canonical continuous algebraic
isomorphism $F(X)\to F_\rho (X)$, and
it can be shown without serious difficulties that the topological group
$F_\rho (X)$ is  complete if and only if
$X$ is Dieudonn\'e complete \cite{Tk3}.
Sipacheva has proved that the free topological group $F(X)$ has a
base of neighbourhoods of identity that are closed in the
topology of the topological group $F_\rho (X)$.

\smallskip

Let $X$ be a set, and let $\Cal V$ and $\Cal W$ be any
two uniformilities on $X$ generating the same Tychonoff topology.
(Such a triple $(X, \Cal V, \Cal W)$ is  termed sometimes
a {\it bi-uniform space}.)

\proclaim{Question}
Does there exist a topological group $F(X,\Cal V, \Cal W)$
algebraically generated by $X$ (free over $X$)
such that $\Cal V$ is the restriction to $X$ of the
{\it left} uniform structure of $G$, and
$\Cal W$ is the restriction to $X$ of the
{\it right} uniform structure?
\endproclaim

This question can be obviously reformulated in terms of
universal arrows to forgetful functors.
This concept may help to understand how the completeness works.

Among other results on the algebro-topological structure
of the free topological groups, let us mention a nice theorem
of Tkachenko \cite{Tk1,2} stating that the free topological group
over a compact space has the c.c.c. property (together with its
subsequent generalization due to Uspenski\u\i\ \cite{U1}), and
a characterization of such Tychonoff spaces $X$ that the
free topological group $F(X)$ embeds into a direct product of a family
of separable metrizable groups \cite{Gur}.

\subheading{
5. Abelian case
}
All of the above results have, of course, their analogs for
free Abelian topological groups. Moreover,
one can also give a very convenient and simple description of topology of
$A(X)$ which has no analog (yet?) in non-commutative case.
One can define Graev metrics on $A(X)$ in the same way as for
$F(X)$, and it turns out that they describe the topology of $A(X)$.
It follows from this observation that the canonical morphism
from $A(X)$ to the free locally convex space $L(X)$ is
an embedding of $A(X)$ as a closed topological subgroup \cite{Tk3}.
Both completeness of $A(X)$ over Dieudonn\'e complete
spaces $X$ and the Abelian analog of the
subgroup theorem were established much earlier than their
non-Abelian counterparts \cite{Tk3}.

The embedding $A(X)\hookrightarrow L(X)$
enables one to describe the topology of $A(X)$ as the topology of
uniform convergence on all equicontinuous families of
characters of $A(X)$, and this way Pontryagin-van Kampen duality
comes into being. For the first time the Pontryagin-van Kampen duality
for free Abelian topological groups was studied by Nickolas \cite{Nic2}
who has shown, answering a question by Noble \cite{No},
that the topological group $A[0,1]$ is non-reflexive
(that is, does not verify the statement of Pontryagin duality theorem).
Later the author had obtained the following  result.

\proclaim{2.9. Theorem {\rm (Pestov \cite{Pe8})}}
Let $X$ be a Dieudonn\'e complete $k$-space with $dim~X=0$.
Then the free Abelian topological group $A(X)$ is reflexive.
\qed\endproclaim

The free topological group $F(X)$ is so ``regularly shaped'' that one
is wondering whether it satisfies any known version of
noncommutative duality (Tannaka-Krein duality is
known to be insufficient) or, at the very least,
whether its topology can be described with the help of
equicontinuous families of homomorphisms from $F(X)$ into
some fixed topological group --- say, $GL(\infty)=
\varinjlim GL(n,\R)$.

\smallpagebreak
\heading   3. M-equivalence and dimension
\endheading
\smallpagebreak

 In 1945 Markov in one of his important
papers \cite{Mar2} asked whether
any two Tychonoff topological spaces, $X$ and $Y$, with isomorphic free
topological groups $F(X)$ and $F(Y)$ are necessarily homeomorphic.
Soon Graev in his no less important papers
\cite{Gr1,2} answered in the negative by constructing a whole series of
pairs $X,Y$ of spaces with $F(X)\cong F(Y)$, therefore the
resulting relation of equivalence between Tychonoff spaces turned
out to be substantial.
Graev called such spaces $X$ and $Y$ {\it F-equivalent}; we follow
the terminology due to Arhangel'ski\u\i\ \cite{Arh3,5,6,8} and call
such spaces {\it Markov equivalent} or {M-equivalent}.
Graev paid special attention to the pairs of spaces
$X,Y$ with Graev free topological groups isomorphic, $F_G(X)\cong F_G(Y)$;
however, the distinction between the two relations of equivalence
is --- from the viewpoint of their topological
properties --- inessential.
With the help of Arhangel'ski\u\i 's terminology, one of the central
results of the Graev's paper \cite{} can be formulated like this.

\proclaim{3.1. Theorem {\rm (1948, Graev)}} If $X$ and $Y$ are $M$-equivalent
compact metrizable spaces then $dim~X=dim~Y$.
\qed\endproclaim

\noindent
(Here $dim~X$ stands for the Lebesgue covering dimension of a space $X$.)

This result --- as well as technique of the proof ---
has received a lot of attention later. The generalizations of the result
came in two directions: firstly, the equivalence relation was
being replaced by  more and more loose ones, and secondly,
the topological restrictions on the spaces $X,Y$ were weakened.

In 1976 Joiner \cite{Joi} noticed that the conclusion $dim~X=dim~Y$
remains true if $X$ and $Y$ are both {\it locally compact metrizable}
spaces such that the free {\it Abelian} topological groups,
$A(X)$ and $A(Y)$, are isomorphic. (Following Arhangel'ski\u\i ,
we call such spaces $X,Y$ {\it A-equivalent}.)
Of course, $A$-equivalence of two topological spaces follows from their
$M$-equivalence, because the universal arrow $A(X)$ is a composition
of the universal arrow $F(X)$ and the functor of abelianization
$\bold{TopGrp}\to\bold{AbTopGrp}$.

Consider the universal arrow from the free
Abelian topological group $A(X)$
to the forgetful functor from the category of locally convex spaces
with weak topology to $\bold{AbTopGrp}$.
The composition of two universal arrows is obviously the
free locally convex space in weak topology, $L_p(X)$.
Therefore, we come to a still looser relation of equivalence between two
spaces: $X$ and $Y$ are {\it l-equivalent} if $L_p(X)\cong L_p(Y)$.
In 1980 Pavlovski\u\i\ \cite{Pa} had shown that $dim~X=dim~Y$ if
$X$ and $Y$ are $l$-{\it equivalent} spaces which are {\it
either locally compact
metrizable} or {\it separable complete metrizable}.

So far all proofs relied on a suitable refinement of the original
Graev's techniques.
A basically new method --- that of inverse spectra --- was invoked
and applied to this problem by Arhangel'ski\u\i\ \cite{Arh3,6} who deduced from
the Pavlovski\u\i 's theorem the following landmark result.

\proclaim{3.2. Theorem {\rm (Arhangel'ski\u\i\ 1980)}}
Let $X$ and $Y$ be $l$-equivalent compact spaces. Then $dim~X=dim~Y$.
\qed\endproclaim

Independently a weaker version was obtained by Zambakhidze \cite{Zam1}:
the covering
dimension of any two $M$-{\it equivalent compact} spaces is the same.
Later this result was generalized by him to the class of
{\it \u Cech complete, scaly, normal, totally paracompact} spaces \cite{Zam2}
(it remained not quite clear how wide this class actually was).
About the same time the result had been
independently somewhat generalized by Valov and
Pasynkov \cite{VP}.

Further efforts have been boosted by
 a question asked by Arhangel'ski\u\i\  \cite{Arh5}:
{\it is it true that for Tychonoff $M$-equivalent
spaces $X$ and $Y$ one has $dim~X=dim~Y$?}

The answer ``yes'' came from the author, who proved in late 1981
\cite{Pe3} the following
result by combining and adjusting both Graev's lemma and the
spectral technique of Arhangel'ski\u\i :

\proclaim{3.3. Theorem {\rm (Pestov, 1981)}}
If $X$ and $Y$ are $l$-equivalent Tychonoff spaces then $dim~X=dim~Y$.
\qed\endproclaim

As a matter of fact, the aforementioned Graev's lemma, which
forms the core of the proofs,
is not a
single result but rather
a {\it scheme} of results, improved and adjusted from one situation
to another.
We present it as it appears
in \cite{Pe5}, not in the most general form possible, but in a quite
elegant one.

\proclaim{3.4. Graev's Lemma}
If $X$ and $Y$ are $M$-equivalent Tychonoff spaces then
$X$ is a union of countably many subspaces each of which is
homeomorphic to a subspace of $Y$.
\qed\endproclaim

Then one is using addition theorems for covering
dimension valid for
spaces with countable base; to proceed from such spaces to
a general situation, the Tychonoff spaces $X$ and $Y$ are decomposed in
inverse spectra of spaces with countable base and the same dimension
as $dim~X$ and $dim~Y$; the property of $l$-equivalence of the
two limit
spaces is partially delegated to the spectrum spaces, in a form
strong enough to
ensure a version of the Graev's lemma.

It was shown by Burov \cite{Bu} that the result and the scheme of the proof
remain true also for cohomological dimension $dim_G$ where the group
of coefficients $G$ is a finitely generated Abelian group
(in particular, $dim_\Z X\equiv dim~X$).

The weak dual space to $L_p(X)$ is the space of continuous functions
on $X$ with the topology of simple (pointwise) convergence, $C_p(X)$
(it follows actually from a version of the Yoneda lemma).
The theory of linear topological structure of the LCS $C_p(X)$
has grown out of Banach space theory, after the following
observation proved useful \cite{Cor}:
any Banach space $E$ in weak topology is a subspace of
$C_p(X)$ where $X$ is the closed unit ball of the dual to $E$
with weak$^\star$ topology.
This theory is developing now on its own, and a good survey
is \cite{Arh9}.
A bridge between theory of spaces $C_p(X)$ and universal arrow theory
is erected by means of the following observation:
since the two LCS's in weak topology,
$L_p(X)$ and $C_p(X)$, are in duality, then
two topological spaces $X$ and $Y$ are $l$-equivalent if and only if
$C_p(X)$ and $C_p(Y)$ are isomorphic.

Arhangel'ski\u\i\  was first to suggest an even weaker realtion
of equivalence between two Tychonoff topological spaces, $X$ and $Y$:
two such spaces are called {\it u-equivalent} if
the locally convex spaces $C_p(X)$ and $C_p(Y)$ are isomorphic as
uniform spaces (with the natural additive uniformity).
Surprisingly, it was possible to make one more step in extending the
original Graev's result.

\proclaim{3.5. Theorem {\rm (Gul'ko, \cite{Gu})}}
If $X$ and $Y$ are $l$-equivalent Tychonoff spaces then $dim~X=dim~Y$.
\qed\endproclaim

The proof of Gul'ko's result \cite{Gu} develops along the same lines as
the author's earlier theorem, but technically it
is considerably more complicated.

One can consider even weaker realtion of equivalence:
two topological spaces, $X$ and $Y$, are said to be
{\it t-equivalent} \cite{GuKh}
if the locally convex spaces $C_p(X)$ and $C_p(Y)$ are {\it homeomorphic}
as topological spaces.
It is not known whether the dimension is
preserved under the relation of $t$-equivalence.
It is worth mentioning that all the aforementioned equivalence
relations
(those of $M$-, $l$-, $u$-, $t$-equivalence)
have been distinguished from each other.

What remains still unclear, is the existence of
a reasonable straightforward
characterization of dimension of a Tychonoff space $X$ in terms
of the additive uniformity of the LCS $C_p(X)$, or  the
linear topological structure of the space $L_p(X)$, or ---
at the very least --- the algebro-topological
strucuture of  $F(X)$.
The existing proofs  are in a sense
obscure and do not reveal the real machinery keeping
dimension preserved by the equivalence relations.

It is an opinion of the author
that emerged from discussions with Gul'ko in 1991 that
a complete understanding of the phenomenon
of preservation of
dimension is to be saught on the following way.

\proclaim{Conjecture}
The Lebesgue
dimension of $X$ can be expressed in terms of a certain
(co)ho\-mo\-lo\-gy theory associated with the LCS in weak topology $L_p(X)$.
\endproclaim

It is not clear if one can use any of the already existing
(co)homology theories for locally convex spaces,
because a desired theory should make a sharp distinction between
week and normable topologies.
For instance, the space $C(X)$ endowed with the topology of
{\it uniform convergence on compacta} instead of the
pointwise topology carries essentially no information about the
dimension of $X$, according to celebrated Milyutin isomorphic
classification theorem \cite{Muly}.

The following remarkable
theorem by Pavlovski\u\i\ may be also sugggestive;
to our knowledge, no attempt has been made yet
to generalize it to arbitrary
$CW$-complexes.

\proclaim{3.6. Theorem {\rm (Pavlovski\u\i\ \cite{Pa})}}
Two polyhedra (simplicial complexes) $X$ and $Y$ are
$l$-equivalent if and only if $dim~X=dim~Y$.
\qed\endproclaim

It is well known that $dim~X\leq n$ if and only if
every continuous mapping from $X$ to a sphere $S^{n+1}$
is homotopically trivial \cite{vM}.
The structure of the free topological groups on spheres
is well understood \cite{KatzMoN2}, so the following sounds sensible.

\proclaim{Conjecture} Let $X$ be a Tychonoff topological space.
Then $dim~X\leq n$ if and only if every continuous
homomorphism from the free topological group $F(X)$ to
the free topological group $F(S^{n+1})$ is homotopically trivial.
\endproclaim

Of course, similar considerations no longer work for $l$-equivalence
because any LCS is contractible,
but the above conjecture may help to reach a deeper understanding
for the relation
of $M$-equivalence.

In addition to Gel'fand-Na\u\i mark duality,
general interest to the problem of preservation of
properties of topological spaces by different
functors from the category $\bold{Tych}$ to the categories of
topological algebra has been heated for a long time by the following
result of Nagata \cite{Nag}.

\proclaim{3.7. Theorem {\rm (Nagata)}}
Two Tychonoff spaces $X$ and $Y$ are homeomorphic if and only if
the topological rings $C_p(X)$ and $C_p(Y)$ are isomorphic.
In other terms, the functor $C_p(\cdot)$ from $\bold{Tych}$
to $\bold{TopRings}$ is a (contravariant) inclusion functor.
\qed\endproclaim

By considering for every Tychonoff space $X$ the universal arrow
from $X$ to a forgetful functor from the category $\bold{TopGrp}$
to $\bold{Tych}$
sending a topological group to a topological subspace
consisting of all elements of order $2$, one comes to the following result
\cite{Pe12}.

\proclaim{3.8. Theorem}
There exists a (covariant) inclusion functor
$\bold{Tych}\to \bold{TopGrp}$.
\qed\endproclaim

There is no {\it full} inclusion functor of such kind \cite{Pe12}.

The following question seems very natural in connection with
our problematics, and it was asked independently by many
(for example, by Zarichny\u\i\ in Baku-1987):

\proclaim{Question}
Is it true that $K$-groups of $M$-equivalent Tychonoff topological spaces
are isomorphic?
\endproclaim

An obvious idea, to obtain the affirmative answer
with the help of universal classifying groups,
fails, because if $G$ is a non-Abelian topological group and
$X$ and $Y$ are $M$-equivalent,
then it follows (from the Yoneda's lemma, actually)
that $K(X)$ and $K(Y)$ are isomorphic as {\it sets},
not {\it groups}: the set $Hom_c(F(X), G)$ does not carry a natural
groups structure because of non-commutativity of $G$ ---
and the universal classifying groups in $K$-theory are noncommutitive.
(This is why a corresponding statement in \cite{VP} is wrong.)

The general classification of topological spaces up to an
$M$-equivalence (as well as $l$-equivelence and other relations
mentioned in this section)
seems a totally hopeless problem.
For numerous results on preservation and non-preservation of
particular properties of set-theoretic topology
by $M$-equivalence, $l$-equivalence etc. see
\cite{Arh3,5,7,8, Gr1,2, Ok, Tch1,2}.
{}From our point of view, there are at least two cases where
such a classification may be achieved.
The first is the case of $l$-equivalence of {\it $CW$-complexes}
 (in view of the
Pavlovski\u\i 's theorem), and the second is the case
of $M$-equivalence of the so-called {\it scattered spaces} \cite{ArhPo}
(in view of the complete classification of
all countable metric spaces up to $M$-equivalence
obtained by Graev \cite{Gr1,2}).

\smallpagebreak
\heading  4. Applications to general topological groups
\endheading
\smallpagebreak

In this section we consider some applications of free
topological groups to general theory of topological groups.
Remark that perhaps one owes the very existence of the concept of
free topological group
to a stimulating applied problem of such kind:
in his historical note \cite{Mar1} Markov was openly guided by the idea
of constructing the first ever example of a
Hausdorff topological group whose
underlying space was not normal.
(The free topological group over any Tychonoff non-normal space $X$
is such.)

Free topological groups provide flexible ``building blocks''
for erecting more sophisticated constructions.
Also, the following theorem is of crucial importance.

\proclaim{4.1. Theorem {\rm (Arhangel'ski\u\i\ \cite{Arh1})}}
Let $f$ be a quotient mapping from a topological space $X$ onto
a topological group $G$. Then the continuous homomorphism
$\hat f: F(X)\to G$ extending $f$ is open and therefore
$G$ is a topological quotient group of $F(X)$.
\qed\endproclaim

Seemingly, analogs of this theorem exist for other types of universal
arrows as well, and one is wondering whether this result can be given
a universal categorial shaping. This result (and its analogs)
are invaluable for examining questions of existence of couniversal objects
of one or another kind.

\subheading{
1. NSS property
}
Our first example is the NSS property. A topological group
$G$ {\it has no small subgroups} if there is a neighbourhood $U$ of the
identity element $e$
such that the only subgroup in $U$ is $\{e\}$.
This is abbreviated to {\it NSS}. The crucial
role of the NSS property in
Lie theory (especially in connection with Hilbert's Fifth Problem) is well
known.

In 1971  Kaplansky wrote (\cite{Kap}, p.89):
{\it ``The following appeaps to be open: if $G$ is NSS and $H$ is a closed
normal subgroup of $G$, is $G/H$ NSS? This is true if in addition
$G$ is locally compact, but we shall only be able to prove it late in the game.
(Of course it is an old result for Lie groups.)''}

Very soon Morris \cite{Mo3} answered in negative by constructing a
counter-example,
and later he and Thompson \cite{MoTh2} have
presented the following

\proclaim{4.2. Theorem}
Let $X$ be a submetrizable Tychonoff topological space (that is,
a Tychonoff space
admitting a continuous metric).  Then the Markov free topological group
$F(X)$ over $X$  is an NSS group.
\qed\endproclaim

It was asked in \cite{MoTh2} whether the following result is true.

\proclaim{4.3. Theorem}
Each topological group is a quotient group
of an NSS group.
\qed\endproclaim

The author \cite{Pe1,4} has deduced Theorem 4.3
from Theorems 4.2 and 4.1
(and later it turned out that such a deduction
follows at once from the above Theorems 4.2 and 4.1
in conjunction with \cite{Ju},
see \cite{Arh4,5}).

It was shown by Sipacheva and Uspenski\u\i\ \cite{SiU}
that both the original proof of Theorem 4.2 by
Morris and Thompson \cite{MoTh2}
and the later proof proposed by Thompson \cite{Th}
are not free of certain deficiencies. In
the same work
\cite{SiU} an elaborate  proof of Theorem 4.2
(definitely ``hard'' --- it relied on combinatorial
technique of words in free groups) was given. Thus,
both results remain valid.
The concept of free
Banach-Lie algebra enables us to provide a purely Lie-theoretic
(and certainly ``soft'') proof of Theorem  4.2
(see Section 7  below).

\subheading{
2. Zero-dimensionality
}
Our next story is about quotient groups of zero-dimensional
topological groups, and it is strikingly similar to the preceding
development. In 1938 Weil (see the note \cite{Arh4} for this and the
next references)
claimed that open continuous homomorphisms of topological groups
do not increase dimension.
This statement was later refuted by Kaplan
by means of a counterexample. Arhangel'ski\u\i\ \cite{Arh2} has shown that
{\it every topological group with countable base is
a quotient group of a zero-dimensional group}.
(Zero-dimensionality here
and in the sequel is understood in the sense of Lebesgue covering
dimension $dim$.) Possible ways to represent {\it any} topological group
as a quotient group of a zero-dimensional one were discussed by
Arhangel'ski\u\i\ in \cite{Arh1}, but it was until late 1980 that the
above conjecture remained open.

\proclaim{4.4. Theorem {\rm (Arhangel'ski\u\i\ \cite{Arh4,5})}}
Any topological group is a topological quotinet group of a group
$G$ with $dim~G=0$.
\qed\endproclaim

Subtle topological considerations involving Graev metrics on free groups
played a crucial role in the proof of the main auxiliary result:
{\it if a submetrizable topological space $X$ is a disjoint union of
a family of spaces each of which has a unique non-isolated point
then $dim~F(X)=0$}. Then the fact that every Tychonoff space
is a quotient of a space with the above property is used,
together with Argangel'ski\u\i 's Theorem 4.1.

This result brought to life a variety of satellite theorems
and examples refining the statement. Of them the most
important one is, from the author's viewpoint, the following.

\proclaim{4.5. Theorem {\rm (Sipacheva \cite{Si2})}}
If $X$ is a Tychonoff space and $dim~X=0$ then $dim~F(X)=0$.
\qed\endproclaim

\subheading{
3. Topologizing a group
}
As the last example, we discuss a problem by Markov \cite{Mar2}
remaining open for 40 years. A subset $X$ of a group $G$ is called
{\it unconditionally closed} in $G$ if $X$ is closed with respect
to every Hausdorff group topology on $G$.
Markov asked \cite{Mar2}
{\it whether a group $G$ admits a connected group topology if
and only if every unconditionally closed subgroup of $G$
has index $\geq\frak c$.}
(Obviously, this condition is necessary.)

The first counterexample was constructed by the author in \cite{Pe13}.
Denote by $L^\flat (X)$ the universal arrow from a uniform space
$X$ to the forgetful functor
from the category of pairs
$(E,Y)$, $E$ a LCS and $Y$ a bounded subset of $E$
(with obviously defined morphisms), to $\bold{Unif}$,
of the form $(E,Y)\to Y$ where $Y$ inherits the
additive uniformily from $E$. If $G$ is a topological group and
$H$ a closed subgroup, then the left action of $G$ on the quotient
space $G/H$ with a natural quotient uniform structure \cite{RoeD}
lifts to a continuous action of $G$ on $L^\flat (G/H)$.
The double semidirect product

$$G^\uparrow=(G\ltimes L^\flat (G) )\ltimes L^\flat (X),$$

\noindent where $X$ is the disjoint sum of a family of copies of
a quotient space of $G\ltimes L^\flat (G)$,
serves as a counterexample to the Markov question
in case where $G$ is an uncountable totally disconnected
topological group.

Later it was observed by Remus \cite{Re} that the
infinite symmetric group $S(X)$ with pointwise topology provides
another --- much more transparent ---
 counterexample to the Markov's conjecture.

The author's techniques was also used by him
to construct an example of a group
admitting a nontrivial Hausdorff group topology but admitting no non-trivial
Hausdorff metrizable topology \cite{Pe11}.

Another problem of Markov still remains open.
A subset $X$ of a topological group $G$ is called
{\it absolutely closed} if it is closed in the coarsest topology
on $G$ making all mappings of the form

$$x\mapsto w(x)$$

\noindent continuous as soon as $w(x)$ is a word in the alphabet formed
by all elements of $G$ and a single variable $x$.
This topology is an analog of the Zariski topology in affine spaces;
we think it is natural to call it the {\it Markov topology} on a group.

\proclaim{Problem {\rm (Markov \cite{Mar2})}}
Prove of refute the conjecture: every
unconditionally closed subset
of a group is absolutely closed.
\endproclaim

Denote by $\frak T_M(G)$ the Markov topology on a group $G$,
and by $\frak T_\wedge (G)$ --- the topology
intersection of all
Hausdorff group topologies on $G$.
It is clear that $\frak T_M(G)\subset \frak T_\wedge (G)$.
The Markov's problem can be be now put in other terms:
{\it is it true that for an arbitrary group $G$ one has
$\frak T_M(G)=\frak T_\wedge (G)$?}

\smallpagebreak
\heading  5. Free products of topological groups
\endheading
\smallpagebreak

Graev \cite{Gr3} presented a constructive description of the topology
of the free product $G\ast H$ of two compact groups;
also he proved a version of Kurosh subgroup theorem in the same paper.
Later both results have been generalized to
$k_\omega$-groups (or, more precisely, topological groups
whose underlying spaces are $k_\omega$) \cite{MoOTh}.
It is known that those results are no longer true beyond the class
of such topological groups.

One can ask about the free products of topological groups almost
the same natural questions as for free topological groups:
to give a reasonable description of topology in general case,
to prove (or refute) that the free product of two (an arbitrary
family of) complete topological groups is a complete group;
to prove (or refute) that if $H_\alpha$ is a topological
subgroup of $G_\alpha$
for every $\alpha\in A$ then $\ast_{\alpha\in A}H_\alpha$ is a topological
subgroup of $\ast_{\alpha\in A}G_\alpha$.
However, here is a question deserving, from our viewpoint,
a special attention --- and not only because of its respectable age.

As it is well known, the construction of free product of groups
is a generalization of the construction of a free group: indeed,
the free group $F(X)$ over the set $X$ of free generators is
just the free product $\ast_{x\in X}\Z_x$ of $\vert X \vert$
copies of the infinite cyclic group $\Z$.
This is obviously not the case with free {\it topological} groups
and free products of {\it topological} groups --- unless
$X$ is discrete.
In 1950 Graev mentioned this and remarked that {\it
``the question of existence of
a natural construction which would embrace both free topological groups
and free products of topological groups still remains open.''}
It does --- for some 42 years already.

Let $\{G_x:x\in X\}$ be a family of topological groups
indexed with elements
of a topological space $X$. One would like to define the free
product $\ast_{x\in X}G_x$ as an appropriate universal arrow
in such a way that
1) in case where $G_x\cong \Z$ for all $x\in X$, the group
$\ast_{x\in X}G_x$ was (naturally) topologically isomorphic to
the Markov free topological group $F(X)$;
2) in case where $X$ is a discrete topological space,
$\ast_{x\in X}G_x$ was a usual free product of topological groups.

Our suggestion is that a clue to the above problem
might be the space $\Cal L(G)$ of all closed
subgroups of a topological group $G$, endowed with an appropriate
topology.
This space (and, moreover, a topological lattice)
has been thoroughly studied \cite{Pr1}
in connection with  extending the Mal'cev
Local Theorems  to the case of locally compact groups.
It is known that there exist numerous ``natural'' topologies
on the set $\Cal L(G)$, including the Vietoris, Chabeauty, and other
topologies ({\it loc. cit.}).

\proclaim{Conjecture}
The free product of a family of topological groups $\{G_x:x\in X\}$
indexed with
elements of a Tychonoff topological space $X$ can be defined as the
universal arrow from $\{G_x:x\in X\}$ to the
functor from the category $\bold{TopGrp}$ to the category of
all families of topological groups indexed with elements of
Tychonoff spaces (with relevant morphisms
between them), which assigns to a topological group $G$ the family
$\{H: H\in \Cal L(G)\}$, the space $\Cal L(G)$ being endowed with
an appropriate topology.
\endproclaim

The Graev problem can be put in connection with deformation theory and
quantum groups.
In quantum physics, one considers {\it deformations}
of algebro-topological objects (such as Lie groups)
as families of objects, $A_{\hbar}$, depending on a continuous
parameter $\hbar$, which is assumed to be a ``very small''
real number approaching zero. Physically, $\hbar$ is the Planck's constant,
and the case $\hbar=0$ corresponds to the (quasi) classical
limit of a theory;
what is deformed, is the object $A_0$.
The absence of nontrivial deformations for classical simple Lie groups and
algebras was a reason for introducing new kind of
objects  --- the quantum groups
\cite{Drin, Man2, RTF, Ros, Wo}.

While there exists a rich mathematically sound
deformation theory for {\it Lie
algebras}, deformations of Lie group are often treated at a
heuristic level. The conjectural Graev construction would
enable one to consider the family $G_{\hbar},~\hbar\geq 0$ of Lie groups
as a veritable
continuous path in the topological space $\Cal L(\ast_{x\in X}G_x)$.

Quantum groups were introduced in mathematical physics
to describe the so-called broken symmetries of physical systems.
The concept of a quantum group is not something accomplished,
and its development is still in progress.
It is only natural, in search of more interrelations between newly explored
categories of mathematical physics, to
look for universal arrows between them.
Does the notion of a free quantum group over a ``quantum space''
make sense?

\smallpagebreak
\heading  7. Free Banach-Lie algebras and their Lie groups
\endheading
\smallpagebreak

The {\it free Banach-Lie} algebra, $\frak{lie}(E)$, over
a normed space $E$ is the universal arrow from $E$ to the forgetful
functor $S$ from the category $\bold{ BLA}$ of complete
Lie algebras endowed with submultiplicative norm to the category $\bold{Norm}$
of normed linear spaces.

\proclaim{7.1. Theorem {\rm (Pestov \cite{Pe18})}}
The free Banach-Lie algebra exists for every normed space $E$,
and $E\hookrightarrow \frak{lie}(E)$ is an isometric embedding.
The Lie algebra $\frak{lie}(E)$ is centerless and infinite-dimensional
if $dim~E>0$.
\qed\endproclaim

One can also define the free Banach-Lie algebra over an arbitrary
pointed metric
space $X$ (we will denote it $\frak{lie}_X$) as the universal arrow from
$X$ to the forgetful functor from $\bold{ BLA}$ to $\bold{Met}_\ast$
(zero goes to the marked point). Obviously, it is just the composition
of the free Banach space and free Banach-Lie algebra arrows.

A Banach-Lie algebra $\frak g$ is called {\it enlargable} if
it comes from a Banach-Lie group.
Every free Banach-Lie algebra is enlargable, and we will denote the
correposnding simply connected Banach-Lie group by
$\frak{LG}(E)$ (resp. $\frak{LG}_X$).
Since every Banach-Lie algebra $\g$ is a quotient Banach-Lie algebra
of the free Banach-Lie algebra over the underlying Banach space
of $\g$, then we come to an independent proof of a
result due to van Est and \'Swierczkowski \cite{\'S3}:
{\it every Banach-Lie algebra is a quotient of an enlargable
Banach-Lie algebra.}

This result can be strengthened. The couniversality of
the Banach space $l_1$ among all separable Banach spaces is well-known
\cite{LiT} (actually, it is due to the fact that
$l_1$ is the free Banach space over a discrete metric space).
Therefore,  $\frak{lie}(l_1)$ is a couniversal
separable Banach-Lie algebra, and the universality
 property is transfered to the
Lie group $\frak{LG}(l_1)$.

\proclaim{7.2. Theorem}
There exists a couniversal connected separable Banach-Lie group.
\qed\endproclaim

Of course, the same is true for groups containing a dense subset of
cardinality $\leq\tau$.

One can show using results of Mycielski \cite{My} and an idea of
Gel'baum \cite{Gel}
that for any metric space $X$, the exponential image of
$X\backslash\{0\}$ in the Lie group $\frak{GL}_X$ generates an
algebraically free subgroup.
Now let $Y$ be a submetrizable pointed space admitting a one-to-one
continuous mapping to $X$. The composition of this mapping and
the
exponential mapping $exp_{\frak{GL}_X}$ determines a
continuous monomorphism $F_G(Y)\to \frak{GL}_X$, and since any
Banach-Lie group has NSS property then it is shared by $F_G(Y)$.
This is the promised ``soft''
proof of Morris-Thompson-Sipacheva-Uspenski\u\i\
theorem.

In view of the existence of a couniversal
separable Banach-Lie group, the following question seems most natural.

\proclaim{Question}
Does there exist a universal separable
Banach-Lie group?
\endproclaim

One should compare it with the following fascinating result
of Uspenski\u\i\ \cite{U4}.

\proclaim{7.3. Theorem {\rm (Uspenski\u\i)}}
The group of isometries of the Banach space $C(I^{\aleph_0})$
endowed with the strong operator topology
 is a universal topological group with countable base.
\qed\endproclaim

However, the general linear group $GL(E)$ of any Banach space
$E$, endowed with the uniform
operator topology,
cannot serve as a universal Banach-Lie group
because there exist enlargable separable Banach-Lie algebras $\g$
which do not admit a faithful linear representation in a Banach space
\cite{vE\'S}.

The universal arrow from a Lie algebra, $\g$, to the forgetful functor
from the category of associative algebras to the category of Lie algebras
is well-known; this is the universal enveloping algebra, $U(\g)$, of $\g$
\cite{Dr}.

It seems that little is known about a topologized version of this, that is,
the universal arrow from a locally convex Lie algebra, $\g$,
to the forgetful functor
from the category of locally convex
associative algebras to the category of locally convex Lie algebras.
Let us denote this arrow by $i_\g:\g\to U_{\Cal T}(\g)$.
Is $i_\g$ an embedding of topological algebras?
(That is, does a topological version of the Poincar\'e-Birkhoff-Witt theorem
hold?)
Is $U_{\Cal T}(\g)$ algebraically isomorphic to $U(\g)$?
What about the convergence of the exponential mapping for $U_{\Cal T}(\g)$?

The only result I am aware of in this connection is
the following.

\proclaim{7.4. Theorem {\rm \cite{Bou}}}
The universal
enveloping algebra $U(\g)$ of a finite-dimensional Lie algebra $\g$ can be
made into a normed algebra if and only if $\g$ is nilpotent.
\qed\endproclaim

This means that, firstly,
a metric version of the universal arrow makes no sense and, secondly,
in general the algebra $U_{\Cal T}(\g)$ is non-normable even if $\g$ is
finite-dimensional.

A detailed analysis of the structure of the locally convex associative algebra
$U_{\Cal T}(\g)$ would be helpful in connection with enlargability problems
for $\g$.

\smallpagebreak
\heading 8.  Lie-Cartan theorem
\endheading
\smallpagebreak

The Lie-Cartan theorem  says that
finite-dimen\-si\-on\-al Lie algebras are enlargable,
and it seems that {\it the question on existence of a ``direct''
proof of the Lie-Cartan theorem,
which would be independent of both known proofs} (the cohomological one by
Cartan \cite{C} and the representation-theoretic one by Ado \cite{Ad}),
{\it is still open.}
For a detailed discussion, see the book \cite{Po}, where
it is claimed that the above question
for a long time received an attention from both
French and Moscow schools of Lie theorists (including Serre).

In this Section we discuss the idea of a conjectural proof
based entirely on universal arrows type constructions
(free topological groups and free Banach-Lie algebras).

It is well known how by means of the
Hausdorff series $H(x,y)$
one can associate in the most natural and straightforward way
a local Lie group (or, rather, a Lie group germ in the sense of \cite{Ro})
to any Banach-Lie algebra $\g$ \cite{Bou}. This is why, according to
a result by \'Swierczkowski \cite{\'S1},
the problem of enlarging a given Banach-Lie algebra $\g$
is completely reduced to the problem of embedding
a local Banach-Lie group $U$
into a topological group $G$ as a
local topological subgroup.

Let $\g$ be a Banach-Lie algebra.
Fix a neighbourhood of zero, $U$, such that the Hausdorff series $H(x,y)$
converges for every $x,y\in U$.
(For example, set $U$ equal to a closed ball of
radius less than $(1/3) log~(3/2)$  \cite{Bou}.)
Denote by $\Cal N_\g$ a closed normal subgroup
generated by all elements of
the form $x^{-1}[x.(-y)]y,~x,y\in U$.
Clearly, the subgroup $\Cal N_\g$ is normal in $F(\g)$
and does not depend on
the particular choice of $U$.
Denote by $G_\g$ the topological group quotient of
$F(\g)$ by $\Cal N_\g$, and
by $\phi_\g:\g\to G_\g$ the restriction of the quotient homomorphism
$\pi_g : F(\g)\to G_\g$ to $\g$.
One can prove that $\pi_g$ is a universal arrow of a certain type.

It is well known (in  different terms, though ---
\cite{\'S2}) that
the enlargability of $\g$ is equivalent to any of the following
conditions:
a) the intersection $\Cal N_\g \cap \g$ is discrete in $\g$;
b) the restriction of $\phi_\g$ to a neighbourhood of zero
in $\g$ is one-to-one;
c) the topological group $G_\g$ can be given a structure of an
analytical Banach-Lie group in such a way that $\phi_\g$ is a local
analytical diffeomorphism; in this case $Lie~(G_\g)\cong\g$,
$\phi_\g = exp_{G_\g}$, and $G_\g$ is simply connected.

Although one can show that
the closedness of $\Cal N_\g$ in general is not sufficent
for any of these conditions to be fulfilled,
it {\it is} so in the following particular case.

\proclaim{8.1. Theorem {\rm \cite{Pe19}}}
A Banach-Lie algebra $\g$ with finite-dimensional
center is enlargable if and only if
the subgroup $\Cal N_\g$ is closed in $F(\g)$.
In this case the quotient topological group $G_\g$ carries a
natural structure of a Banach-Lie group associated to $\g$.
\qed\endproclaim

The proof of this result goes as follows: firstly, it is reduced to
separable Banach-Lie algebras with the help of a local theorem
\cite{Pe14}, and then certain perfectly
direct and functorial constructions are used, including
the free Banach-Lie algebra over the underlying Banach space of
$\g$, the  Banach-Lie group associated to it, and their quotients.

Now only one obstacle remains between us and a direct proof of the
Lie-Cartan theorem.

\proclaim{Conjecture}
The closedness of the subgroup
$\Cal N_\g$ in the free topological group
$F(\g)$ over the underlying topological space
of a finite-dimensional Lie algebra $\g$ can be proved
relying solely on the description of topology
of free topological groups over finite-dimensional Euclidean spaces.
\endproclaim

The subgroup $\Cal N_\g$ is compactly generated; since the
compact set
generating it is in $F_3(\g)$ rather than $F_1(\g)$ then
one should single out some additional algebro-topological
property of the group $\Cal N_\g$ which would ensure the
closedness (or completeness).

We already {\it know} that $\Cal N_\g$ {\it is always closed in} $F(\g)$
for $\g$ finite-dimensional (it follows from the Lie-Cartan theorem),
and the problem looks so natural in this setting.
It is so tempting to think that
the genuine reason why the statement of
Lie-Cartan theorem is always true for finite
dimensional Lie algebras, is not (co)homological but
entirely in the realm of general topology, namely:
finite dimensional Lie algebras are $k_\omega$ spaces,
while infinite dimensional ones are not.

\medpagebreak
\heading  9. Locally convex Lie algebras and groups
\endheading
\smallpagebreak

Infinite-dimensional groups play a major role in the contemporary
pure and applied mathematics \cite{Kac1,2}.
Many of them cannot be given a structure of a Banach-Lie group
(for example, groups of diffeomorphisms of manifolds,
some of their subgroups preserving a certain differential-geometric
structure,
Kac-Moody groups).
At the same time, in all particular examples
to an infinite-dimensional group there is associated in some natural
way an infinite-dimensional Lie algebra, and therefore
it is appealing, to try to develop a version of Lie theory
with all its attributes general enough to embrace all
particular examples of infinite-dimensional groups.

Such attempts have lead to the theory of Lie groups modeled over
locally convex spaces (bornological and sequentially complete \cite{Mil}),
especially over Lie groups modeled over Fr\'echet spaces \cite{KoYMO}.
We will call by a {\it Fr\'echet-Lie group}
a group object in the category of smooth Fr\'echet manifolds,
that is --- in this case --- just
a smooth manifold modeled
over a Fr\'echet space which carries a group structure such that the
group operations are Fr\'echet $C^\infty$.

There is a striking difference between the Banach and Fr\'echet versions of
Lie theory. For example, although there is a well-defined notion of
the Lie algebra, $Lie~(G)$,
 of a Fr\'echet-Lie group $G$ (which is a
Fr\'echet-Lie algebra), the exponential mapping $exp_G: Lie~(G)\to G$
need not be $C^\infty$
nor a local diffeomorphism;  therefore there is in general no canonical
atlas on a Fr\'echet-Lie group.
Moreover, the following question seems to be still open:

\proclaim{Question {\rm \cite{Mil, KoYMO}}}
Does the exponential map $exp_G: Lie~(G)\to G$
always exist for a Fr\'echet-Lie group $G$?
\endproclaim

Because of such misbehaviour of Fr\'echet-Lie theory,
some mathematicians are
questioning its ability to serve
as a basis for
infinite-dimensional group theory.
Among them is Kirillov who once
(Novosibirsk, January 1988) even expressed the
opinion that obtaining an answer to
the above question either in positive or in
negative sense would be disadvantageous all the same!

Nevertheless,
we believe that this question should be answered
in order to understand the proper place of
Fr\'echet-Lie theory, and now
we want to present a new, universal arrow type, construction of
locally convex Lie algebras, which may give a clue.

It is convenient to present the results in the spirit of
$\Delta$-normed spaces and algebras belonging to
Antonovski\u\i , Boltyanski\u\i\ and Sarymsakov \cite{ABS}.

Let $\Delta$ be a directed partially ordered set.
A vector space $E$ is said to be $\Delta$-normed if
there is fixed a family of seminorms
$p=\{p_\delta :\delta\in\Delta\}$ with the property
$p_\delta\leq p_\gamma \Leftrightarrow \delta\leq\gamma$.
(The family $p$
is called a $\Delta${\it -norm} because it
can be treated as a single map $E\times E\to \R^\Delta$
where $\R^\Delta$ is the so-called {\it topological semifield},
and it satisfies close analogs of all three axioms of a usual
norm.)

Let $A$ be an algebra. We will say that a
$\Delta$-norm $p=\{p_\delta :\delta\in\Delta\}$ on $A$
is {\it submultiplicative} if

\item{(i)} {\it for every $\delta,\gamma\in\Delta$ such that
$\delta<\gamma$ and for every $x,y\in A$ one has
$p_\delta (x\ast y)\leq p_\gamma (x)\cdot p_\gamma (y)$,
where $\ast$ denotes the binary algebra operation};

\item{(ii)} {\it for every $\delta\in\Delta$
there is a $\gamma$ such that for every $x,y\in A$ one has
$p_\delta (x\ast y)\leq p_\gamma (x)\cdot p_\gamma (y)$.}

One can show that the topology of
every locally convex topological algebra
is given by an appropriate submultiplicative $\Delta$-norm.
For example, the locally multiplicatively convex
topological algebras introduced by Arens and Michael \cite{Ar, Mic}
are characterized by existence of a $\Delta$-norm with the
property $p_\delta (x\ast y)\leq p_\delta (x)\cdot p_\delta (y)$
for all $x,y\in A$ and every $\delta\in\Delta$.

For a fixed directed set $\Delta$ the class of all complete
$\Delta$-normed Lie algebras forms a category with
contracting Lie algebra homomorphisms as morphisms.
We will denote this category $\Delta\bold{LA}$.

\proclaim{9.1. Theorem}
For every $\Delta$-normed vector space $(E,p)$ there exists a
universal arrow from this space to the forgetful functor
from $\Delta\bold{LA}$ to the category of $\Delta$-normed spaces.
It is an isometric embedding of $(E,p)$ into a
$\Delta$-submultiplicateively normed Lie algebra $\frak{lie}(E)$.
\qed\endproclaim

In a particular case where $\Delta$ is a one-point set,
the above construction coincides with the construction of a
free Banach-Lie algebra over a normed space considered earlier.

If $\Delta$ has countable cofinality type
(in particular, is countable) then the Lie algebra $\frak{lie}(E)$
is a Fr\'echet-Lie algebra.

The algebra $\frak{lie}(E)$ is centerless and infinite-dimensional
(unless $dim~E=1$).
It is completely unclear whether such Fr\'echet-Lie algebras are
enlargable (that is, come from Fr\'echet-Lie groups).
The property of being centerless gives a hope that the answer is
``yes,'' at least in some cases.
However, if $\Delta=\N$ and
a corresponding sequence of seminorms, $p$,
 grows ``fast enough,'' there is a good
evidence that $\frak{lie} (E,p)$  can have no exponential map.

\proclaim{9.2. Theorem}
Let $(E,\norm\cdot )$ be a normed space.
Define a $\Delta$-norm $p$, where $\Delta=\N$, by
letting $p_n=n!\norm\cdot,~n\in\N$.
Suppose there exist a Fr\'echet-Lie group, $G$, associated to
the Lie algebra $\frak{lie}(E)$. Then there is no exponential map
$\frak{lie}(E)\to G$.
\qed\endproclaim

One can also study {\it free locally convex Lie algebras} over
locally convex spaces, that is, universal arrows from an
LCS $E$ to the forgetful functor from the category of locally convex
topological Lie algebras and continuous Lie algebra homomorphisms
to the category of locally convex spaces.
We will denote the free localy convex Lie algebra over $E$ by
$\frak{LClie}(E)$.
If $X$ is a Tychonoff space, then one can consider
the {\it free locally convex Lie algebra} over $X$, defined either
as the coomposition of the free locally convex space $L(X)$
and the free locally convex Lie algebra, or
directly as the universal arrow from $X$ to
the forgetful functor from the category of locally convex
topological Lie algebras and continuous Lie algebra homomorphisms
to the category $\bold{Tych}$. We denote this Lie algebra
by $\frak{LClie}_X$.

P. de la Harpe has kindly drawn my attention to the following problematics.

\proclaim{Problem {\rm (Bourbaki \cite{Bou})}}
Is it true that every extension of a Lie algebra $\g$
by means of a $\g$-module $M$ is trivial
(in other terms, $H^2(\g,M)=(0)$ for
every $\g$-module $M$)
if and only if
$\g$ is a free Lie algebra?
\endproclaim

The property $H^2(\g,M)=(0)$ is readily verifiable for a free Lie algebra
$\g$, but the validity of inverse implication is not known.

It is not clear yet whether free locally convex Lie algebras can help
in answering the above question (supposedly in negative), but
at the very least,
they enjoy a similar property for {\it continuous} second cohomology.

\proclaim{9.3. Theorem}
Let $X$ be a separable metrizable topological space, and let
$M$ be a complete normable locally convex $\frak{LClie}_X$-module.
Then every locally convex extension of the Lie algebra
$\frak{LClie}_X$ by means of $M$ is trivial. In particular,
$H_c^2(\frak{LClie}_X,M)=(0)$.
\qed\endproclaim

The proof follows the argument for free Lie algebras,
but the Michael Selection Theorem (Theorem 1.4.9 in \cite{vM})
is involved.

In some cases one managed to establish the triviality of
{\it algebraic} second cohomology for
locally convex (and even Banach) Lie algebras \cite{dlH}.

\smallpagebreak
\heading  10. Supermathematics
\endheading
\smallpagebreak

The (unhappy but hardly avoidable)
term ``supermathematics'' is used to designate the mathematical
background of dynamical theories with nontrivial fermionic sector
in the quasi-classical limit
$\hbar\to 0$.
The ``supermathematics'' includes
superalgebra, superanalysis, supergeometry etc., all of these being
obtained from their ``ordinary'' counterparts by incorporating
odd (anticommuting) quantities
\cite{BBHR, BBHRPe, B, DeW, Man1}.

In one of those approaches an important role is played by the so-called
ground algebras, or algebras of supernumbers; in other approach,
algebras of this type come into being as algebras of superfunctions over
purely odd supermanifolds. As a matter of fact, those algebras
turn out to be universal arrows of a special kind, and they also
find an independent application in infinite-dimensional
differential geometry.

We will give necessary definitions.
The term {\it ``graded''} in this paper means {\it ``${\Bbb Z}_2$-graded''}.
A graded algebra $\Lambda$ is
an associative algebra over the
basic field $\Bbb K$ together
with a fixed vector space decomposition
$\Lambda \cong \Lambda^0 \oplus \Lambda^1$,
where $\Lambda^0$ is called the {\it even} and $\Lambda^1$ the {\it odd
part (sector)} of $\Lambda$, in such a way that the {\it parity}
$\tilde x$ of any element $x \in \Lambda^0 \cup \Lambda^1$
defined by letting $x \in \Lambda^{\tilde x}, \tilde x \in
\lbrace 0,1 \rbrace = \Bbb Z_2$, meets the following restriction:

$$\tilde{xy} = \tilde x + \tilde y,~~  x,y \in \Lambda^0 \cup
\Lambda^1  $$

If in addition one has

$$ xy = (-1)^{\tilde x\tilde y}yx,~~ x,y \in \Lambda^0 \cup
\Lambda^1  $$

then $\Lambda$ is called {\it graded commutative}.

\proclaim{10.1. Theorem {\rm \cite{Pe16,17}}}
Let $E$ be a normed space. There exists a universal arrow
$\wedge_B E$ from $E$ to the forgetful functor from the category
of complete submultiplicatively normed graded commutative algebras
to the category of normed spaces.
It contains $B$ as a
normed subspace of the odd part $(\wedge_B E)^1$
in such a way that $E \cap \{1\}$ topologically generates
$\wedge_BE$ and every linear operator $f$ from $E$ to the odd part $\Lambda^1$
of a complete normed associative unital graded commutative algebra
$\Lambda$ with a norm $\Vert f \Vert_{op} \leq 1$ extends to an
even homomorphism $\hat f \: \wedge_B E \to \Lambda$ with a norm
$\Vert \hat f \Vert_{op} \leq 1$.
\qed\endproclaim

Algebraically, $\wedge_B E$ is just
the exterior algebra over the space $E$,
endowed with a relevant norm and completed after that.
It enjoys one more property.
A {\it Banach-Grassmann algebra} \cite{JP} is a complete normed
associative unital graded commutative algebra $\Lambda$
satisfying the following two conditions.

$BG_1$ {\it (Jadczyk-Pilch self-duality).}
For any $r,s \in \Bbb Z_2 = \{0,1\}$ and any bounded
$\Lambda^0$-linear operator $T\: \Lambda^r \to \Lambda^s$
there exists a unique element $a \in \Lambda^{r+s}$
such that $Tx = ax$ whenever $x \in \Lambda^r$.
In addition, $\Vert a \Vert$ equals the operator norm $\Vert T \Vert_{op}$
of $T$.

$BG_2$. The algebra $\Lambda$ decomposes into an $l_1$ type sum
$\Lambda \simeq \Bbb K \oplus J^0_{\Lambda} \oplus \Lambda^1$
where $\Bbb K = \Bbb R$ or $\Bbb C$ and $J^0_{\Lambda}$ is
the even part of the closed ideal $J_{\Lambda}$
topologically generated by the odd part $\Lambda^1$.
In other words, for an arbitrary $x \in \Lambda$ there exist elements
$x_B \in \Bbb K, ~ x_S^0 \in J^0_{\Lambda}$, and
$x^1 \in \Lambda^1$ such that $x = x_B + x^0_S + x^1$ and
$\Vert x \Vert = \Vert x_B \Vert + \Vert x^0_S \Vert +
\Vert x^1 \Vert$.

\proclaim{10.2. Theorem {\rm \cite{Pe17}}}
Let $E$ be a normed space.
The following conditions are equivalent:

\item{(i)} $dim~E=0$;
\item{(ii)} $\wedge_B E$ is a Banach-Grassmann algebra.
\qed\endproclaim

The algebra $\wedge_B l_1$ (denoted by $B_\infty$)
was widely used in superanalysis
\cite{JP}.

The algebras of the type $\wedge_B E$ appear in infinite-dimensional
differential geometry: in \cite{KL}, Klimek and Lesniewski used them
for constructing Pfaffian systems over infinite-dimensional Banach spaces
after it became clear that the earlier considered Pfaffians over
Hilbert spaces are insufficient for applications in mathematical
physics.

If one wishes to study algebras of superfunctions on
purely odd (that is, including fermionic degrees
of freedom only) infinite dimensional supermanifolds modeled
over {\it locally convex} spaces, then another
universal arrow comes into being.
A {\it locally convex} graded algebra $\Lambda$
 carries two structures - that of
a graded algebra and of locally convex space --- in such a way that
multiplication is continuous and both even and odd sectors are closed
subspaces of $\Lambda$.
A topological algebra $A$ is called
{\it locally multiplicatively convex}, or just
{\it locally m-convex}, if
its topology can be described by a family of all
submultiplicative continuous seminorms. (Equivalently:
$A$ can be embedded into the direct product of family of
normable topological algebra.) \cite{Ar, Mic}
An {\it Arens-Michael algebra}  \cite{He}
is a complete locally m-convex algebra.

\proclaim{10.3. Theorem {\rm \cite{Pe15,16}}}
Let $E$ be a locally convex space. Then
there exists a universal arrow
$\wedge_{AM}E$ from $E$ to the forgetful functor from the category
of graded commutative Arens-Michael algebras to the category of
locally convex spaces.
\qed\endproclaim

The two particular cases are well-known:
$\wedge_{AM}\R^{\aleph_0}$ is the DeWitt supernumber algebra \cite{DeW},
and $\wedge_{AM}\R^{\omega}$ is the nuclear $(LB)$ algebra
considered in \cite{KoN}.
(Here $\R^{\aleph_0}$ stands for the direct product of
countably many copies of $\R$, and $\R^{\omega}$ denotes the direct limit
$\varinjlim\R^n$.) In addition, in the finite-dimensional case,
$\wedge_{AM}\R^q$ is just the Grassmann algebra with $q$
odd generators.

Perhaps, the same sort of construction would serve as a base for
study of Pfafians on infinite dimensional locally convex spaces.

At present one of the most appealing unsolved problem in
``supermathematics'' is
to give a unified treatment
of all existing approaches to the notion of a supermanifold
by viewing supermanifolds over non-trivial ground algebras
$\Lambda$ as superbundles over $Spec~\Lambda$.

Denote by $\Cal G$  the category of finite-dimensional Grassmann algebras
and unital algebra homomorphisms preserving the grading.
Let $\bold{LCS}^{{\Cal G}^{op}}$
denote the category of all contravariant functors from
$\Cal G$ to the category  $\bold{LCS}$ of locally convex spaces and
continuous linear operators;
 we call the category
$\bold{LCS}^{{\Cal G}^{op}}$ the category of {\it virtual
locally convex superspaces}.
Every graded locally convex space $E=E^0\oplus E^1$ canonically
becomes an object of  $\bold{LCS}^{{\Cal G}^{op}}$,
because it determines a functor of the form
$\wedge (q) \mapsto [\wedge (q)\otimes E]^0$;
we will identify this functor with
$E$. The simplest nontrivial
example of a virtual graded locally convex space is $\R^{1,1}=\R^1\oplus \R^1$.
The category $\bold{LCS}^{{\Cal G}^{op}}$ is a subcategory of the category
$\bold{DiffLCS}^{{\Cal G}^{op}}$ of
all contravariant functors from
$\Cal G$ to the category  $\bold{DiffLCS}$ of locally convex spaces and
infinitely smooth mappings between them.

\proclaim{Conjecture}
The set of all morphisms in the category $\bold{DiffLCS}^{{\Cal G}^{op}}$
from a purely odd graded
locally convex space $E$ to $\R^{1,1}$ carries a natural structure of
a graded locally convex algebra canonically isomorphic to the
free graded commutative Arens-Michael algebra, $\wedge_{AM}E'_\beta$, on
the strong dual space $E'_\beta$.
\endproclaim

\smallpagebreak
\heading
11. $C^\star$ algebras and noncommutative mathematics
\endheading
\smallpagebreak

Every normed space $E$ admits a universal arrow to the forgetful
functor from the category of (commutative)
$C^\star$-algebras and their morphisms
to the category of normed spaces and
contracting linear operators; we will denote it by $C^\star (E)$
($C_{com}^\star (E)$, in commutative case),
and refer to as {\it the free (commutative) $C^\star$-algebra
over a normed space}.
The arrows in both cases are isomorphic embeddings.
This is simply due to the two facts:
firstly, every normed space $E$ embeds
into the $C^\star$-algebra of continuous functions on
the closed unit ball of the dual space $E'$ with the weak$^\star$ topology,
and secondly, the class of (commutative) $C^\star$ algebras is
closed under the $l_\infty$-type sum.

This construction is a particular case of the
Blackadar's construction of a $C^\star$-algebra
defined by generators and relations \cite{Bla}.
For example, the {\it free $C^\star$-algebra}
over a set $\Gamma$ of free generators \cite{GM}
is just the free $C^\star$-algebra in our sense over the Banach space
$l_1(\Gamma)$.
In non-commutative topology \cite{BOB}
the $C^\star$-algebras $C^\star (l_1(\Gamma))$
(treated as objects of the opposite category)
are viewed as
noncommutative versions of Tychonoff cubes $I^\tau$,
because they are couniversal objects
(universal --- in the opposite category).

It is known that every free $C^\star$-algebra is {\it residually
finite-dimensional (RFD)},
that is, admits a family of $C^\star$-algebra homomorphisms
to finite-dimensional $C^\star$-algebras separating points \cite{GM}.
The same is true for our more general objects.

\proclaim{11.1. Theorem}
For every normed space $E$ the $C^\star$-algebra
$C^\star (E)$ is residually finite-dimensional.
\qed\endproclaim

This result seems interesting because there are few known classes of
RFD $C^\star$-algebras \cite{ExL}.

Both embeddings have been considered earlier \cite{BlP, Ru}, where
the so-called matrix norms on $E$ defined by those embeddings
are denoted by $MAX$ and $MIN$.
This construction is especially important for
the so-called  quantized functional analysis \cite{Eff},
of which the idea is that all the main functional-analytic properties and
results concerning Banach spaces can be expressed in terms of the
universal arrow $C_{com}^\star (E)$,  so
their non-commutative versions stated for
$C^\star (E)$
constitute the object of {\it quantized} (that is, noncommutative)
functional analysis.

In this connection, it may be useful to
consider two equivalence relations on Banach spaces,
two such spaces, $E$ and $F$, being equivalent iff
$C^\star (E)\cong C^\star (F)$
(respectively, $C_{com}^\star (E)\cong C_{com}^\star (F)$).

If one wishes to study ``quantized'' theory of LCS's then
one should turn to the similar universal arrows from a given
LCS $E$ to the forgetful functor from the category of
the so-called {\it pro-$C^\star$-algebras} in the sense
of N.C. Phillips \cite{Ph} (just inverse limits of $C^\star$-algebras)
and their morphisms to the forgetful functor to the category of
LCS's; there are both commutative and non-commutative versions of
those universal arrows.

Finally, we expect that a whole new class of examples of the
so-called
{\it quantum algebras} in the sense of Jaffe and collaborators
\cite{JLO}
can be obtained by considering universal arrows
from a set of data including
graded normed spaces to the relevant forgetful functor.

\smallpagebreak
\heading
Acknowledgments
\endheading
\smallpagebreak

It is an appropriate occasion
to name all Institutions that supported my
investigations related to universal arrows during more than 13 years;
I am grateful to all of them. Those are:
Institute of Applied Mathematics and Mechanics at Tomsk State
University (1978-80) and Faculty of Mechanics and Mathematics
of the same University (1983-91);
Division of Mathematics of Moscow State University (1980-83);
Faculty of Mathematics of Far Eastern State University, Vladivostok (1984);
Institute of Mathematics, Novosibirsk
Science Centre (1988-91);
Department of Mathematics, University of Genoa, and
Group of Mathematical Physics of the
Italian National Research Council (1990 and 1991);
Department of Mathematics, University of Victoria, Canada (1991-92);
and Department of Mathematics, Victoria University of Wellington,
New Zealand (1992---). A 1992 research grant V212/451/RGNT/594/153
from the Internal Grant Committee of the latter University
is acknowledged.

In the present  area of research I was most
influenced by ideas of my Ph.D. advisor Professor A.V. Arhangel'ski\u\i .
My warm thanks also go to  S.M. Berger, U. Bruzzo,
M.M. Choban, R.I. Goldblatt, S.P. Gul'ko,
P. de la Harpe, A.E. Hurd,
S.S. Kutateladze, T. Loring, J. Leslie, N.C. Phillips, M.G. Tkachenko,
O.V. Sipacheva, V.V. Uspenski\u\i\
for their help expressed in various forms.

The present article is an extended version of my invited talk
at the University of Wollongong (Australia),
and I am most grateful to Professor S.A. Morris
for providing such an excellent opportunity.

\vfill
\eject

\bye